\documentclass[11pt]{article}
\usepackage{amssymb, amsmath, amsthm, amsfonts,xcolor, enumerate}
\usepackage{fullpage}
\usepackage{hyperref}

\usepackage{t1enc}
\usepackage[T1]{fontenc}
\usepackage[utf8]{inputenc}

\newtheorem{theorem}{Theorem}[section]
\newtheorem{lemma}[theorem]{Lemma}

\newtheorem{claim}[theorem]{Claim}
\newtheorem{prop}[theorem]{Proposition}
\newcommand{\ep}{\varepsilon}
\newcommand{\cB}{\mathcal{B}}

\def\qedf{\hfill $\Box$}

\def\komment#1{}
\let\komment=\footnote

\title{Regular decomposition of the edge set of graphs with applications}
\author{ B\'ela Csaba\thanks{Bolyai Institute, University of Szeged, Szeged, Hungary. This research was partially supported  by 
NKFIH grant KH 129597 and by the grant NKFIH-1279-2/2020 of the Ministry
for Innovation and Technology, Hungary.} }

\date{}

\begin{document}
\maketitle

\begin{abstract}
We introduce a new method for decomposing the edge set of a graph, and use it to replace the Regularity lemma of Szemer\'edi in some 
graph embedding problems. An algorithmic version is also given. 
\end{abstract}

\section{Introduction}

The Szemerédi Regularity lemma is among the most powerful tools in graph theory. The lemma appeared in~\cite{Szemeredi1} in 1978. A weaker version had already been used earlier 
by Szemer\'edi to prove the Erd\H os-Tur\'an conjecture~\cite{ET, Szemeredi2}. Since then the lemma and its ramifications have found several applications in graph and hypergraph theory,
number theory, algebra, geometry, and computer science. 

  
Given a simple\footnote{We consider only simple graphs in this paper.} graph $G=(V, E)$ and a number $\varepsilon \in (0,1),$ the lemma asserts that $V$
can be partitioned into $k\le N(\ep)$ subsets $V_1 \cup \ldots \cup V_k$ and another set $V_0=V - (\cup_iV_i)$ such that $G[V_i, V_j]$ is $\ep$-regular for every $1\le i \neq j \le k,$
except at most $\ep k^2$ pairs of indices, $|V_0|\le \ep n$ and $|V_i|=(n-|V_0|)/k$ for every $1\le i \le k.$

In the proof the threshold $N(\ep)$ for the number of parts is a tower of twos of height $O(\varepsilon^{-5}).$ As Gowers proved in~\cite{Gowers}, this tower-type bound is 
unavoidable in general, more
precisely, there are graphs for which $N(\ep)$ has to be at least a tower of twos of height $\Omega(\varepsilon^{-1/16}).$ Conlon and Fox~\cite{CF} further improved the lower bound
to $\Omega(\varepsilon^{-1}).$ This shows the major drawback of the Regularity lemma: in order to get meaningful results, one has to work with enormously large graphs, thereby ruling out practical applications of the lemma. 
We remark that if the number of edges between vertex sets $A$ and $B$ is $o(|A|\cdot |B|),$ then the $(A, B)$-pair is $\varepsilon$-regular by definition, unless 
$\varepsilon$ is very small compared to $|A|$ and $|B|.$ 
Hence, the lemma is useful only when the density of the graph is essentially a positive constant.  

In order to avoid these disadvantages, several alternatives for substituting the Regularity lemma have been discovered, e.g., the weak regularity lemma of Frieze and Kannan~\cite{FK}, or the cylindrical regularity lemma by Eaton and Rödl~\cite{ADLRY, E, ER}. 
Let us also mention a regularity lemma by Gowers~\cite{GowersB}, and a somewhat similar one by the author~\cite{Cs1}. There are also versions for dense subgraphs of random graphs or $C_4$-free graphs~\cite{Fox2}, for graphs with bounded VC-dimension~\cite{Alon1, LovSzeg} or for graphs with low threshold-rank~\cite{GT}. This list is long, but still far from being complete. However, none of the above results have the full strength of the original Regularity lemma.   

In this paper we introduce a new, iterative graph decomposition method. Our technique allows one to use some of the methods developed in the applications of the Regularity lemma of Szemer\'edi. 
As we will see, this decomposition may be used for solving embedding or packing problems in graphs on $n$ vertices that have $\Omega(n^2/(\log n)^{1/q})$ edges, where $q> 2$ is a 
number depending on the problem itself. Besides, the number of vertices in the graph is allowed to be just ``reasonably'' 
large. 

In the decomposition of this paper one does not partition the vertex set of the graph as in the Szemer\'edi Regularity lemma, rather the edge set. 
The basic building blocks of the decomposition are pseudorandom subgraphs, usually super-regular pairs or regular pairs. 
The most important difference between this method and the many versions of the Regularity lemma is that the regular pairs of our decomposition are edge-disjoint, but their vertex sets may intersect. These intersections make the use of the edge decomposition method more difficult, 
but it still can replace the Regularity lemma in many cases. We will demonstrate this via some examples. Still, we think that one cannot use this technique for so diverse problems as the Regularity lemma. 

Let us remark that the present decomposition method is similar to the cylindrical regularity lemma for bipartite graphs by Eaton and R\"odl~\cite{E, ER} and~\cite{ADLRY}, although 
it is slightly less restrictive.   
Our main motivation is the graph functional method of Koml\'os~\cite{KS}. This allows us to get better bounds than that of the cylindrical regularity lemma, which
can be especially important in case the graph has vanishing density. 

We also provide an algorithmic version: a decomposition that can be found in polynomial time by a randomized algorithm. 
The algorithmic version works for somewhat larger and somewhat denser graphs than the one based on graph functionals, but it  is still meaningful for graphs having vanishing densities. 

The outline of the paper is as follows. In Section 2  we review the necessary notions, results for bundle decompositions, and state our main theorem. In Section 3 we prove our main theorem via several  lemmas. 
In Section 4 we apply our method for packing large trees and bounded degree bipartite graphs into a host graph such that only $o(n^2)$ edges remain uncovered. In Section 5 we prove a conditional triangle removal lemma for graphs having relatively few $C_5$s. 
Finally, in Section 6 we prove an algorithmic version of our main theorem.

We made no attempts to optimize on the constants in the paper, 
and will not be concerned with floor signs and divisibility in the proofs. This makes the notation simpler, easier to follow. Throughout the paper $\log x$ will denote the {\it natural logarithm} of $x,$ and $exp(f)=e^f$ for any expression $f.$

\section{Preparations}

\subsection{Notations, definitions}\label{}

Given a graph $G=(V, E)$ we use the notations $v(G)=|V|$ and $e(G)=|E|.$ For disjoint subsets $X, Y\subset V$ we let $G[X, Y]$ denote the bipartite subgraph of 
$G$ with parts $X$ and $Y$ that contains all the edges of $G$ with one endpoint in $X$ and the other endpoint in $Y.$ For every vertex $v\in V$ the neighborhood of $v$ is denoted by 
$N(v),$ and the degree of $v$ is denoted by $deg(v)=|N(v)|.$ Given a set $S\subset V$ we let $N(v, S)=N(v)\cap S$ and $deg(v, S)=|N(v, S)|.$
 
The {\it density} of $G$ is defined to be $d_G=e(G)\cdot {\binom{v(G)}{2}}^{-1}.$ The {\it bipartite} density of bipartite subgraphs of $G$ with parts $A$ and $B$ 
is defined to be $d_G(A, B)=\frac{e(G[A, B])}{|A|\cdot |B|}.$ Sometimes the subscript may be omitted when there is no confusion. Similarly, when a graph in question is bipartite, density will mean bipartite density, unless stated otherwise.
We omit the proof of the following well-known result.

\begin{claim}\label{suru}
Let  $G=(V, E)$ be an $n$-vertex graph with density $d_G.$ Then there exists $X\subset V,$ $|X|=\lfloor n/2 \rfloor$ such that $d_G(X, V-X)\ge d_G.$  
\end{claim}

Another folklore result, which is sometimes called {\it convexity of density} (see e.g.~in~\cite{KS}), will also prove to be useful later.

\begin{claim}\label{cnvx} Let $F=F(A, B; E)$ be a bipartite graph, and let $k, m$ be integers such that $1\le k \le |A|$ and $1\le m \le |B|.$ Then 
$$d_F(A, B)=\frac{1}{\binom{|A|}{k}\binom{|B|}{m}}\sum_{X\in \binom{A}{k}, Y\in \binom{B}{m}} d_F(X, Y).$$
\end{claim}

\begin{claim}\label{convexity}
Let $G=(A, B; E)$ be a bipartite graph with parts $A$ and $B.$ Let $k, m$ be integers such that $1\le k \le |A|$ and $1\le m\le |B|.$ Then there exist 
$X\subset A$ with $|X|=m$ and $Y\subset B$ with $|Y|=k$ such that
$d_G(A, B)\le d_G(X, Y).$ 
\end{claim}

\noindent {\bf Proof:} The claim easily seen to follow from Claim~\ref{cnvx}.
\qedf 

\medskip

We say that a bipartite graph $H=(A, B; E)$ is $\varepsilon$-regular for a real number $\varepsilon>0,$ if $$|d(A, B)-d(X, Y)|\le \varepsilon$$ whenever $X\subset A$ and
$Y\subset B$ such that $|X|\ge \varepsilon |A|$ and $|Y|\ge \varepsilon |B|.$ For $\varepsilon, \delta>0$ we call $H$ an $(\varepsilon, \delta)$-{\it super-regular}\footnote{We remark that in some earlier papers super-regularity meant a somewhat different notion: $\varepsilon$-regularity was replaced by the condition that between any two sets $X\subset A,$ $Y\subset B$ with $|X|\ge \varepsilon |A|$ and $|Y|\ge \varepsilon |B|$ we have $d(X, Y)\ge \delta.$ This was also used by Fox~\cite{Fox1} for proving his breakthrough result on the removal lemma. In order to avoid confusion, we call this notion
$(\varepsilon, \delta)$-super-lower-regularity. } pair, 
if it is $\varepsilon$-regular, and every $v\in A$ has at least $\delta |B|$ neighbors and every $u\in B$ has at least $\delta |A|$ neighbors. 

Let $\cB$ denote the class of balanced bipartite graphs, that is, bipartite graphs having equal-sized parts, and let ${\cB}_m$ denote the class of balanced bipartite graphs having $m$ vertices in both parts. 
The claim below follows from the definition of $\ep$-regularity. 

\begin{claim}\label{szup-atlagfok}
Let $0<\ep<1/4,$ and assume that $F\in \cB_m$ is an $\ep$-regular pair with density $d_F=d.$ 

\noindent $(i)$ There exists $H\subset F,$ $H\in \cB_{m'}$
such that $H$ is a $(2\ep, d-2\ep)$-super-regular pair having density $d_H\ge d-2\ep,$ where $m'\ge (1-\ep)m.$

\noindent $(ii)$ Denote the parts of $F$ by $A$ and $B,$ and let $S_A\subset A,$ $S_B\subset B$ such that $|S_A|, |S_B|\ge \ep m +k$ (here $k$ is a positive integer).
Then $S_A$ has at least $k$ vertices $v_1, \ldots, v_k$ such that $deg_F(v_i, S_A)\ge (d-\ep)|S_B|.$ 
\end{claim}

\noindent {\bf Proof:} For proving $(i)$ observe that by $\ep$-regularity, $F$ may have at most $\ep m$ vertices in both its vertex parts which have degree smaller than $(d-\ep)m.$ Those vertices
will be discarded. Then we make the parts having equal sizes by discarding some more (at most $\ep m$) arbitrarily chosen vertices from the larger part. Call the remaining subgraph $H.$
One can easily verify that the statement holds for $H.$

We use the definition of $\ep$-regularity for $(ii)$ as well. Let $S_A'\subset S_A$ be an arbitrary subset with $|S_A'|=\ep m.$ Since $S_A'$ and $S_B$ are large, the density of edges between them is at least $d-\ep.$ Hence, we can find a vertex $v_1\in S_A'$ such that
$v_1$ has at least $(d-\ep)|S_B|$ neighbors in $S_B.$  Next add an arbitrary vertex of $S_A-S_A'$ to $S_A',$ and then delete $v_1$ from $S_A'.$ The above argument can be repeated to find 
further $v_i$ vertices of $S_A$ that have large degree into $S_B.$   
\qedf

\subsection{The edge-decomposition theorem}

We are ready to state our main result:

\begin{theorem}\label{dekomp}
Let $G$ be a balanced bipartite graph on $n$ vertices with density $d_G,$ and let $0<d \le 1$ and $0< \ep <1/12$ such that $n> exp(10 \log (1/d) \log (1/\ep)/\ep^2).$ 
Then the \emph{edge set} of $G$ can be decomposed as follows: 
$E(G)$ can be written as the edge-disjoint union of the $(\varepsilon, d-\varepsilon)$-super-regular 
bipartite graphs $H_1, \ldots, H_{K}\in \cB,$ and another bipartite graph $H_0\in \cB,$ where $K=K(\varepsilon, d).$ For $i\ge 1$ each $H_i$ has at least $m=m(\varepsilon, d)$ vertices in both parts and density at least $d,$ 
while $H_0$ has density less than $d.$ Furthermore, $$m \ge d^{(10/\varepsilon^2)\log (1/\varepsilon)}n/2$$ and 
$$K \le 2d_G \cdot d^{-(20/\varepsilon^2)\log (1/\varepsilon)}/d.$$ 
\end{theorem}


\medskip

Note that while implicit, the density of $G$ plays an important role: if $d_G$ is less than $d,$ then $H_0=G.$ 

The theorem is about balanced bipartite graphs, however, by Claim~\ref{suru} one may use it for arbitrary graphs so that there is no loss in the density. 
Hence, if one deletes the condition that $G$ is a balanced bipartite graph, the resulting statement holds without any further
assumptions. 

As we indicated in the introduction, we prove Theorem~\ref{dekomp} using the graph functional method of Koml\'os. 

\section{Proof of Theorem~\ref{dekomp}}

We begin with the definition of our graph functional. After showing some of its properties, we will prove the theorem. 

\subsection{A graph functional for finding super-regular pairs}

The method of graph functionals was introduced by Koml\'os~\cite{KS}. 
The general form of such a graph functional 
$\Phi$ for some graph $H$ is as follows: $$\Phi(H)=\varphi(d_H)\cdot v(H),$$ where (as before) $d_H$ denotes the density of $H$ and $\varphi$ is an increasing real function.  
We will be interested in the special case $\varphi(x)=x^r/2$ for $r=r(\varepsilon),$ with $0< \varepsilon \ll 1.$ 

\medskip

Given a graph $G=(V, E)$ we
consider the following maximization problem: $$\max_{H\subset G, \ H\in \cB} \Phi(H)=\max_{H\subset G, \ H\in \cB}d_H^r \cdot v(H).$$ We will prove that if $d_G$ is not too
small 
and $r=r(\varepsilon)$ is sufficiently large, then the subgraph $H\subset G$ at which $\Phi$ attains its maximum is a 
super-regular pair. 
 
\begin{prop}\label{en}
Let $G=(V, E)$ be a balanced bipartite graph on $n$ vertices with density $d_G$ and $0<\varepsilon<1/4$ be a real number. Set 
$r= (10/\varepsilon^2)\log (1/\varepsilon).$ If $n> exp(r \log (1/d_G))$ and $d_G>\ep,$ then $G$ contains an 
$(\varepsilon,\delta)$-super-regular subgraph $H\in \cB$
with $\delta \ge d_G-\varepsilon,$ $v(H)\ge d^r_Gn/2$ and density $d_H\ge d_G.$
\end{prop}

Proposition~\ref{en} is proved via a few lemmas. For the first one let's assume that $H=(A, B; E)$ is a balanced bipartite graph such that $m=|A|=|B|.$ Let $\Psi$ be a
graph functional defined on balanced bipartite subgraphs $F$ of $H$ as follows: $$\Psi(F)=d_F^q v(F),$$ where $q\in [2/\ep, m]$ is an arbitrary real number. Notice,
that here we allow the exponent of the density to be smaller than the $r$ in Proposition~\ref{en}.

\begin{lemma}\label{Fok} Assume that there exists $v\in A$ such that $deg(v)\le (d_H-\varepsilon)m.$ Then there is a vertex $w\in B$ such that $\Psi(H')>\Psi(H),$ where $H'$ is
the subgraph of $H$ induced by the vertex classes $A-v$ and $B-w.$
\end{lemma}

\noindent {\bf Proof:} First we take an arbitrary vertex $w\in B$ such that $deg(w)\le d_H m.$ Such a vertex must exist using that $e(H)=d_Hm^2.$ 
Let us compute the density of $H'=H[A-v, B-w].$ Since $H'$ has $m-1$ vertices in both vertex classes, we have
$$d_{H'}\ge\frac{d_Hm^2-2d_Hm+\varepsilon m}{(m-1)^2}.$$ Hence, 
$$\frac{\Psi(H')}{\Psi(H)}\ge \left(\frac{(d_Hm^2-2d_Hm+\varepsilon m)m^2}{(m-1)^2d_Hm^2} \right)^q\frac{m-1}{m}.$$
After simplifying the above expression we get that $$\frac{\Psi(H')}{\Psi(H)} \ge \left( \frac{m^2-2m+\frac{\varepsilon m}{d_H}}{(m-1)^2}  \right)^q \frac{m-1}{m}.$$ We show that the right hand side of the above inequality is larger than 1. This implies that $\Psi(H')>\Psi(H),$ thereby proving the lemma.
This is equivalent to the following
$$\left( \frac{m^2-2m+\frac{\varepsilon m}{d_H}}{(m-1)^2}  \right)^q >1 +\frac{1}{m-1}.$$ Using Bernoulli's inequality\footnote{That is, $(1+h)^n \ge 1+nh,$ if $h\ge -1.$} we have
$$\left( \frac{m^2-2m+\frac{\varepsilon m}{d_H}}{(m-1)^2}  \right)^q = \left( \frac{(m-1)^2 -1 +\frac{\varepsilon m}{d_H}}{(m-1)^2}  \right)^q  \ge 1+q\frac{\frac{\varepsilon m}{d_H}-1}{(m-1)^2}.$$ Using that
$q\in [2/\varepsilon, a]$ and that we must have $0\le d_H\le 1$ one can easily verify that the left hand side is larger than 1, proving what was desired.\qedf


\medskip

From now on we assume that $\Phi$ attains its maximum at $H\subset G$ for some $H\in \cB_m.$ By Lemma~\ref{Fok} we have that every vertex has at least 
$(d_H-\ep)m$ neighbors in $H.$ Let the vertex classes of $H$ be denoted by $A$ and $B.$ 
We will prove that $H$ is an $\ep$-regular pair, together with Lemma~\ref{Fok} this will imply Proposition~\ref{en}. 

Let $A_1\subset A$ and $B_1\subset B$ be arbitrary subsets such that 
$|A_1|, |B_1|= \varepsilon m,$ and let $A_2=A-A_1$ and $B_2=B-B_1.$ 
We will prove that $$-\ep \le d(A_1, B_1)-d_H \le \ep^2/5.$$ Using Claim~\ref{convexity} this implies that 
$H$ is an $\ep$-regular pair.

\begin{lemma}\label{l1} 
Using the above notation we have $d(A_1, B_1) \le d_H(1+\varepsilon^2/5).$ 
\end{lemma}

\noindent {\bf Proof:} Assume on the contrary that $d(A_1, B_1) > d_H(1+\varepsilon^2/5).$ Then 
by definition, $$\Phi(H[A_1, B_1])> \left(d_H(1+\varepsilon^2/5)\right)^r\varepsilon m.$$ 
Since $\Phi(H)=d_H^r m,$ we get that 
$$\frac{\Phi(H[A_1, B_1])}{\Phi(H)}> (1+\varepsilon^2/5)^r\varepsilon.$$ 

Clearly, if we can show that  $(1+\varepsilon^2/5)^r\varepsilon>1,$ we arrive at a contradiction, as $\Phi(H)$ has maximum value.

\begin{claim} We have $(1+\varepsilon^2/5)^r> \frac{1}{\varepsilon}.$ 
\end{claim}

\noindent {\bf Proof of the claim:} Let us take the natural logarithm of both sides of the inequality. Using that $r=10\log (1/\ep)/\ep^2,$ we have to prove the following
$$\frac{10\log(1/\ep) \log(1+\ep^2/5)}{\ep^2} > \log(1/\ep).$$ Since $0<\ep<1/4,$ we can simplify the above, and obtain the inequality
$$10\log(1+\ep^2/5) > \ep^2.$$ For proving this, we consider the function $D(x)=10\log(1+x^2/5)-x^2.$ We have $$D'(x)=2x \left( \frac{2}{1+x^2/5} - 1\right).$$ This
derivative is positive for every $x\in (0,1),$ and since $D(0)=0,$ we get that $D(x)>0$ for every $x\in (0,1),$ proving the claim.\qedf  

\noindent With the above we have got the desired contradiction, so the lemma holds.\qedf

\medskip

Next we prove that the density between $A_2$ and $B_2$ cannot deviate much from the density of $H.$

\begin{lemma}\label{l2}
The density $d(A_2, B_2)<d_H \cdot (1+\varepsilon^3/3).$
\end{lemma}

\noindent {\bf Proof:}
Assume on the contrary that $d(A_2, B_2)\ge d_H (1+\varepsilon^3/3).$ Then 
$$\Phi(H[A_2, B_2])\ge d_H^r(1+\varepsilon^3/3)^r(1-\varepsilon)m.$$ We show that  $$d_H^r(1+\varepsilon^3/3)^r(1-\varepsilon)m>d_H^rm.$$
This is equivalent to the inequality $$(1+\varepsilon^3/3)^r>1/(1-\varepsilon).$$
Since  $(1-\ep)(1+2\ep)=1+\ep -\ep^2$ and $\varepsilon \in (0, 1/4)$ we have $$1/(1-\varepsilon)<1+2\varepsilon.$$ 
Using that $1+x\le e^x,$ we get $$1/(1-\varepsilon)< e^{2\varepsilon}.$$ 
We will show that $$(1+\varepsilon^3/3)^r >e^{2\ep},$$ which implies the lemma.
We need a simple claim:

\begin{claim}
We have $e^{\ep^3/4}< 1+\ep^3/3.$
\end{claim}

\noindent {\bf Proof of the claim:} Using the power series expansion of $e^x$ we get
$$e^{\ep^3/4}=1+\frac{\ep^3}{4}+\left(\frac{\ep^3}{4}\right)^2\frac{1}{2}+ \left(\frac{\ep^3}{4}\right)^3\frac{1}{6}+ \cdots < 1+\frac{\ep^3}{4} +\frac{\ep^6}{2}\left(1+\frac{1}{2}+\frac{1}{4}+ \cdots\right) =1+\frac{\ep^3}{4} + \ep^6.$$ Since $\ep \in (0, 1/4)$ we have that $\ep^6< \ep^3/12,$ implying the claim. 
\qedf

\medskip

Using the above claim and substituting the value of $r$ we obtain that $$(1+\varepsilon^3/3)^r > e^{r \ep^3/4}=e^{2.5 \ep\log(1/\ep)}>e^{2\ep},$$ what was desired.
\qedf

\medskip

We are ready to prove Proposition~\ref{en}.

\smallskip

\noindent {\bf Proof:}
Using Lemmas~\ref{convexity},~\ref{l1} and~\ref{l2} we have that the densities $d_H(A_1, B_2)$ and $d_H(A_2, B_1)$ are less than $d_H(1+\varepsilon^2/5)$ and $d_H(A_2, B_2)\le d_H(1+\varepsilon^3/3).$ Assuming that $d_H(A_1, B_1)\le d_H-\ep,$ we get the following upper bound for the number of edges in $H$:
$$e(H)\le (d_H-\varepsilon)\varepsilon^2m^2+2d_H(1+\varepsilon^2/5)\varepsilon(1-\varepsilon)m^2+d_H(1+\varepsilon^3/3)(1-\varepsilon)^2m^2=$$
$$d_H\left(1+\frac{11}{15}\varepsilon^3 -\frac{16}{15}\varepsilon^4 +\varepsilon^5/3\right)m^2 - \varepsilon^3m^2<d_Hm^2,$$ 
here we used that $\varepsilon \le 1/4.$ This is a contradiction, hence, the density $d(A_1, B_1)$ must be at least $d_H-\varepsilon,$
implying that $H$ is an $(\varepsilon, d_H-\ep)$-super-regular pair.

In order to prove the lower bound for the number of vertices consider the following inequalities:
$$m\ge d_H^r m=\Phi(H)\ge \Phi(G)= d_G^r \cdot n/2,$$ here the first inequality follows from the fact that every density is at most 1, while the second is implied by the maximality of $\Phi(H).$ Finally, the density of $H$ must be at least $d_G,$ since in the above inequality $m\le n/2$ always holds. This finishes the proof of Proposition~\ref{en}.
\qedf

\bigskip

\subsection{Finishing the proof of Theorem~\ref{dekomp}}

Let $G=(V, E)$ be a balanced bipartite graph on $n$ vertices with density $d_G.$ Let $0<\varepsilon<1/4$ and $0<d\le 1$ be real numbers. For finding the decomposition of $G$ we apply Proposition~\ref{en} repeatedly.

In the first step we check if $d_G\ge d.$ If not, we let $H_0=G,$ and stop. Otherwise let $H_1$ be a balanced bipartite graph for which the graph functional 
$\Phi(H)=d_H^r\cdot v(H)$ is maximal. By Proposition~\ref{en} $H_1$ is an $(\ep, d_H-\ep)$-super-regular pair. 
Let $G_1=(V, E_1)$ be a subgraph of $G$ what we obtain by deleting the edges of $H_1$ from $G,$ that is, $E_1=E(G)-E(H_1).$
Next we repeat the above procedure for $G_1.$ First check if the density of $G_1$ is at least $d.$ If not, we let $H_0=G_1.$ If $G_1$ is sufficiently large, then 
we apply Proposition~\ref{en} again in order to obtain a new super-regular pair in $G_1,$ which we call $H_2.$

In general, let us assume that we have found the first $i$ super-regular pairs $H_1, H_2, \ldots, H_{i}.$ Let $G_{i}=(V, E_{i}),$ where $E_i=E(G_{i-1})-E(H_i).$ 
We check if 
$d_{G_i}\ge d.$ If not, then we let $H_0=G_i,$ and stop. Otherwise we use Proposition~\ref{en} for finding $H_{i+1}.$ 
 
Next we consider the bounds for $m$ and $k.$ The lower bound for $m$ is just the lower bound for $v(H)$ in Proposition~\ref{en}, with $d$ plugged in for the density of the graph. 
For the upper bound for $K$ notice, that $Kdm^2\le d_Gn^2,$ since the $H_i$'s are edge-disjoint, each has at least $m$ vertices in both parts and density at least $d,$ while $e(G)=d_G n^2/2.$ Reordering the inequality and using that $m \ge d^{(10/\varepsilon^2)\log (1/\varepsilon)}n/2$ gives the claimed bound: 
$K\le 2d_G \cdot d^{-(20/\varepsilon^2)\log (1/\varepsilon)}/d.$ \hfill\qedf

\section{Applying Proposition~\ref{en} for approximate decompositions}

The study of packing of graphs dates back to more than a century, and recently it has received much 
attention~\cite{Cs2, KuhnBl, Ober}. The generic 
packing question is as follows: Given a ``large'' host graph $G$ and a ``small'' graph $H,$ is it possible to cover the edge set of
$G$ edge-disjointly by copies of $H$? Equivalently, we sometimes say that the edge set of $G$ is decomposed by edge-disjoint copies of $H$. In many cases the single small graph $H$ is replaced by a family of graphs.
Perhaps one of the most beautiful questions in the area is the still open Gyárfás--Lehel conjecture from 1978 on the decomposition of the edge set of
$K_n$ into trees having every order between $1$ and $n.$ 

In general even an approximate form of these problems, when we may leave a small percentage of the edge set of the host graph uncovered, is very challenging. Below we present an  application of Proposition~\ref{en}
for approximately decomposing the edge set of a graph into edge-disjoint trees. Depending on $\ep$ and $\delta,$ these trees may have linear size and linear maximum degree.

\smallskip

Let $T$ be a rooted tree on $t$ vertices, its root is denoted by $r.$ Given any $x\in V(T)$ let $N^*(x)$ denote the set of its children, and $deg_T^*(x)=|N^*(x)|.$ We define the $L_i$ {\it level sets} of $T$ as follows:
$L_1=\{r\},$ and for $i\ge 1$ we have that $L_{i+1}=\cup_{x\in L_i}N^*(x).$ Hence, if a vertex $y$ lies in $L_j,$ then $y$ is at distance $j-1$ from $r.$ Let us denote the total number of levels by $s.$

\begin{lemma}\label{EgyFa}
Let $T$ be as above. Assume that $H(A, B; E)$ is an $(\ep, \delta)$-super-regular pair with vertex parts $A$ and $B,$ such that $m=|A|=|B|\ge 2t,$ and $\delta\ge 3\ep.$ 
We further assume that $|L_i|\le \delta m/4$ for every $i\ge 1.$ Then $T\subset H.$
\end{lemma}

\noindent {\bf Proof:}
Our goal is to find an edge-preserving injection $\varphi: V(T) \longrightarrow A\cup B.$ The algorithm we use constructs $\varphi$ step-by-step, alternately embedding consecutive levels of $T,$ 
starting at $L_1=\{r\}.$ 

At any point in time we denote the uncovered subset of $A$ by $A_u,$ and similarly, the uncovered subset of $B$ by $B_u.$ We need to define two more subsets:
$$A'=\{v\in A_u: \ deg_H(v, B_u)\ge (\delta -\ep)|B_u|\} \quad\mbox{and}\quad B'=\{v\in B_u: \ deg_H(v, A_u)\ge (\delta -\ep)|A_u|\}.$$

Observe, that $A_u, B_u, A'$ and $B'$ shrink dynamically as we embed more and more levels. After succesfully embedding, say, level $L_i$ into $A'$ we update these
sets as follows: $A_u=A_u-\varphi(L_i),$ $A'=A'-\varphi(L_i),$ and we also have to leave out those vertices of $B'$ that have degree less than $(\delta -\ep)|A_u|$ in the newly updated $A_u,$ the set $B_u$ remains the same. We do the updating analogously when a level is embedded into $B'.$

The proof of the claim below follows easily from our discussion above and the definition of super-regular pairs, we leave it for the reader.

\begin{claim}\label{meret}
In the beginning we have $A_u=A'=A,$ $B_u=B'=B.$ 
\end{claim}

\begin{claim}\label{szintek}
Given the embedding of $L_i$ into $A',$ we can embed $L_{i+1}$ into $B'$ greedily, here $1\le i\le s-1.$ Analogous statement holds in case $L_i$ is embedded into $B'.$
\end{claim}

\noindent {\bf Proof:} Assume that $x\in L_i,$ and therefore $v=\varphi(x)\in A'.$ By the definition of $A'$ we have that $deg_G(v, B_u)\ge (\delta-\ep)|B_u| \ge |L_{i+1}|+\ep m.$
Since $|B_u|-|B'|\le \ep m$ using $\ep$-regularity, this implies that $deg_G(v, B')\ge |L_{i+1}|.$ 

Hence, for every $x\in L_i$ we can greedily choose a subset $S_x\subset N(v, B')$ such that $|S_x|=deg^*_T(x),$ and $S_x\cap S_{x'}=\emptyset$ for every $x'\neq x,$ $x'\in L_i.$
Therefore we may choose the  $\varphi(y)$ images for $y\in N^*(x)$ greedily from $S_x.$ This extension of $\varphi$ clearly satisfies the requirements of the lemma. It is easy to
see that the same procedure works if $L_i$ is embedded into $B',$ only the letters $A$ and $B$ has to be switched.  
\qedf

\smallskip

We have the following algorithm for finding the embedding $\varphi.$ First we embed $L_1$ into $A',$ which is possible by Claim~\ref{meret}. Then we repeatedly apply Claim~\ref{szintek}, and embed more and more levels, always using the $B'$ or $A'$ subsets. Since $m\ge 2t,$ we always have sufficiently many vertices so that $|A'|, |B'| \ge \delta m/2 - \ep m.$ Hence, this procedure never gets stuck, proving what was desired. 
\qedf

\medskip

\begin{theorem}\label{deko}
Let $0< \ep < 1/4$ and $3\ep \le \delta\le 1$ be real numbers, $m$ and $k$ be positive integers such that $k\le \delta m/4.$ Set $r=(10/\ep^2)\log 1/\ep.$ Assume that $T_1, T_2, \ldots, T_l$ are rooted trees, each on at 
most $t\le m/2$ vertices. Denote their total number of edges by $e_T=\sum_i e(T_i).$ 
Let $G$ be a balanced bipartite graph on $2n$ vertices with density $d_G$ such that $d_Gn^2\ge e_T+(\delta+\ep)n^2.$ Assume further that every level set of each tree has at most $k$ vertices.
If $n\ge m \cdot exp(r \log 1/(\delta+\ep)),$ then $G$ has edge-disjoint copies of $T_1, \ldots, T_l.$
\end{theorem}

\noindent {\bf Proof:} The theorem follows from Proposition~\ref{en} and Lemma~\ref{EgyFa}. First, we apply Proposition~\ref{en} in order to find an $(\ep, \delta)$-super-regular pair $H_1,$
where $H_1$ has at least $m$ vertices in both of its parts. 
Since $H_1$ and $T_1$ satisfies the requirements of Lemma~\ref{EgyFa}, we can embed a copy of $T_1$ into $H_1.$ We used up precisely $e(T_1)$ edges of $G$ this way, we delete them from $H_1,$
the rest of $H_1$ is added back to the graph. This procedure is repeated for every $i\ge 2:$ we find a large $(\ep, \delta)$-super-regular pair $H_i$ by Proposition~\ref{en} having at least $m$
vertices in both parts, since the density of what is left from $G$ will always be at least $\delta+\ep,$ no matter how many trees were embedded thus far. Then we embed $T_i$ into $H_i$ by Lemma~\ref{EgyFa}, delete
the edges of $H_i$ we used for $T_i,$ and finally, the rest of $H_i$ is added back to the graph. 
\qedf

\medskip

Note that, as we indicated earlier, the maximum degree of a tree can be even linear, if $\ep, \delta$ are (small) constants. 

The idea of the above almost decomposition, that is, iteratively applying Proposition~\ref{en} for the vacant part of the host graph after embedding a subgraph, can be used for edge-disjoint packing of other kind of graphs. 

\begin{theorem}\label{deko2} Let $G$ be a graph on $n$ vertices with density $d=\Omega(1).$ Assume that $\ep, \delta$ are real numbers and $D$ a positive integer such that $0<\ep < \delta^{5D}<\delta \ll d.$ 
Set $r=10\log (1/\ep)/\ep^2.$ Let $H_1, \ldots, H_{\ell}$ be bipartite graphs with 
maximum degree $\Delta(H_i)\le D$ for every $i,$ and assume that $v(H_i)\le d^rn/2.$ If $\sum_i e(H_i)\le (1-2\delta)e(G),$ then $H_1, \ldots, H_{\ell}$ can be packed edge-disjointly into $G.$
\end{theorem}

For the proof we need the Blow-up lemma~\cite{KSSz} of Komlós, Sárközy and Szemerédi.

\begin{theorem}\label{Blowup}[Blow-up Lemma] 
Given a graph $R$ of order $r$ and positive integers $\delta$ and $\Delta$
there exists a positive $\ep=(\delta, \Delta, r)$ such that the following holds: Let $n_1, n_2, \ldots, n_r$ be
arbitrary positive parameters and let us replace the vertices $v_1, v_2, \ldots, v_r$ of $R$ with pairwise
disjoint sets $W_1, W_2, \ldots, W_r$ of sizes $n_1, n_2, \ldots, n_r$ (blowing up $R$). We construct two graphs
on the same vertex set $V = \cup_i W_i$: The first graph $F$ is obtained by replacing each edge
$v_iv_j \in  E(R)$ with the complete bipartite graph between $W_i$ and $W_j.$ A sparser graph $G$ is
constructed by replacing each edge $v_iv_j$ arbitrarily with an $(\ep,\delta)$-super-regular pair between
$W_i$ and $W_j.$ If a graph $H$ with maximum degree $\Delta(H)\le \Delta$  is embeddable into $F$ then it is already embeddable
into $G.$
\end{theorem}

\smallskip

We are going to use the Blow-up lemma for the case when $R$ contains two vertices ($r=2$) and the edge that connects them. For us it is more useful to have a more explicit relation 
of $\ep$ and $\delta$ than what is given in the above formulation. By carefully reading the proof of the Blow-up lemma in~\cite{KSSz},
one can see that for the $r=2$ case it is sufficient if $\ep \le \delta^{5\Delta}\le 1.$ 

\medskip

\noindent {\bf Proof: } Assume that we have already packed the first $i$ graphs $H_1, H_2, \ldots, H_i.$ The density of what is left from $G$ must still be at least $\delta.$ Using Proposition~\ref{en} one can find a large $(\ep, \delta)$-super-regular pair in $G.$ Using the Blow-up lemma one can find a copy of $H_{i+1}.$ If $i+1=\ell,$ then we are done.
If not, the remaining density in the host graph $G$ is large enough to repeat the procedure.
\qedf 

\medskip

Theorem~\ref{deko2} gives an approximate solution of a special case of the Oberwolfach problem (see e.g.~in~\cite{Ober}, where the solution is also given), in which 2-regular spanning 
graphs are packed into $K_n.$ While we only pack subgraphs that are relatively small compared to the host graph (but still can have linear size), our host graph does not have to be the complete graph, 
and the subgraphs to be packed can be arbitrary bounded degree bipartite graphs.    


\section{A conditional triangle removal lemma}

The celebrated result of Ruzsa and Szemer\'edi~\cite{RSz} states, roughly speaking, that if a graph of order $n$ has $o(n^3)$ triangles, then it can be made triangle-free
by deleting $o(n^2)$ edges from it. 
Let us give a more precise formulation:

\begin{theorem}
Let $G$ be a graph on $n$ vertices. Then for every $\ep>0$ there exists $\delta=\delta(\ep) > 0$ such that if $G$ has  at most $\delta n^3$ triangles, then
it can be made triangle-free by removing at most $\ep n^2$ edges.
\end{theorem}

This very important result has far reaching implications in graph and hypergraph theory, number theory, etc., see for example in~\cite{CF2}. The currently best bound for $\delta$
is a tower of twos of height $O(\log {1/\ep})$ by Fox~\cite{Fox1}.

Bollobás and Győri raised the following question in~\cite{BGy}: How many triangles can a graph $G$ on $n$ vertices have, if $G$ has no cycle of length 5? They proved that the number of triangles in such a graph is at most $(5/4)n^{3/2}+o(n^{3/2}).$ Since then several improvements were found on the constant multiplier of $n^{3/2}$~\cite{Gy, EM}.

In~\cite{Fox2} Conlon et al.~considered a conditional removal problem, which relates the above two questions. They proved that, given a graph of order $n$ with $o(n^2)$ copies of $C_4$'s and
$o(n^{5/2})$ copies of $C_5$'s, it can be made $\{C_3, C_5\}$-free by deleting $o(n^{3/2})$ edges. Below we state and prove a somewhat similar theorem in which we have 
no condition on the number of $C_4$'s of the graph, only for the number of $C_5$'s. Since $C_4$'s emerge when the number of edges is $\Omega(n^{3/2}),$
the bound for the number of edges to be removed is much larger. Still, the $\delta$ below is just a single exponential function of a polynomial of $1/\ep,$ a huge gain compared 
to the $\delta$ in the known proofs of the (unconditional) triangle removal lemma.

\begin{theorem}\label{removal}
Let $G$ be a tripartite graph with vertex parts $X, Y$ and $Z$ such that $|X|=|Y|=|Z|=n.$ Let $0<\ep \le 1/4$ and set 
$K_0=exp((20/\varepsilon^2)\log^2 (1/\varepsilon))/\ep$ and $m=exp(-10\log^2 (1/\varepsilon)/\varepsilon^2)n.$ If the number of $C_5$'s in $G$ is at most 
$$\delta n^5=\frac{\ep^6 m^2n^3}{2K_0^2},$$ then $G$ can be made triangle-free by deleting at most $4\ep n^2$ edges. 
\end{theorem}

\noindent{\bf Proof:} We begin with applying Theorem~\ref{dekomp} with parameters $\ep$ and $d=2\ep$ for the bipartite subgraph $G[X, Y].$ 
We obtain the $(\ep, \ep)$-super-regular pairs $H_1, H_2, \ldots, H_K,$
each having parts on at least $m$ vertices, and $H_0$ with less than $dn^2=2\ep n^2$ edges. Note that $K\le K_0.$

Call a $C_5=x_1y_1x_2y_2z$ {\it good}, if $x_1, x_2\in X,$ $y_1, y_2\in Y,$ $z\in Z,$  its edges are $x_1y_1, y_1x_2, x_2y_2, y_2z, zx_1,$ and $x_1y_1, y_1x_2, x_2y_2$ belong to the same super-regular pair of the
decomposition. Observe that if $x'y'z'$ is a triangle in $G-H_0$ such that $x'\in X, y'\in Y$ and $z'\in Z,$ then the $x'y'$ edge must belong to a bundle.

Simple computation tells that the number of $z\in Z$ which appears in more than $$\frac{\ep^5 m^2n^2}{2K_0^2}=\frac{1}{\ep}\frac{\ep^6 m^2n^2}{2K_0^2}$$ good $C_5$'s is less than $\ep n.$ 
Call such vertices of $Z$ {\it bad,} the rest of $Z$ are the {\it good} vertices. Delete every edge that joins a bad vertex to any vertex in $X,$ this way we deleted
at most $\ep n^2$ edges. 

Let $z\in Z$ be any good vertex. Assume that it has at least $\ep n/K$ neighbors in both parts of some $H_i.$ Denote the parts of $H_i$ by $X_i\subset X$ and
$Y_i\subset Y,$ and let $m_i=|X_i|=|Y_i|\ge m.$ Notice that if $x\in X_i,$ then by the 
$(\ep, \ep)$-super-regularity it has at least $\ep m_i$ neighbors in $Y,$ and similar holds for any $y\in Y_i.$ Hence, if $x\in X_i, y\in Y_i$ and $x, y\in N(x),$
then, using $\ep$-regularity, there are at least $\ep^3m_i^2$ path of length 3 that connects $x$ and $y$ in $H_i.$ Using our assumption that $z$ has at least $\ep n/K$ neighbors in $X_i$ and $Y_i,$ there are at least $\frac{\ep^5 m^2n^2}{K^2}$ such good $C_5$'s which contain $z$ and three edges from $H_i.$ Since this number
is larger than the upper bound of good $C_5$'s which may contain $z,$ we conclude that for every good vertex $z\in Z$ there is no $H_i$ ($1 \le i\le K$) for which $z$ has many
neighbors in both parts of it. 

Next we repeat the following for every good $z\in Z,$ for every $1\le i \le K$: if $deg(z, X_i)\le deg(z, Y_i),$ then delete all edges that join $z$ to $X_i,$ otherwise delete 
all the edges that join $z$ to $Y_i.$ This way we delete at most $n\cdot K\cdot \ep n/K=\ep n^2$ edges from $G.$ Since after these deletions no good $C_5$ is left, we have removed every triangle. 
\qedf 

\medskip

We remark that in~\cite{Fox2} it was proved that there exist $n$-vertex graphs with $o(n^{2.442})$ $C_5$’s that cannot be made triangle-free
by deleting $o(n^{3/2})$ edges. Note the large gap between the two bounds for conditional removal.

The proof method of Theorem~\ref{removal} can easily be generalized to prove statements of the following type. Let $2\le k\le \ell$ be integers. If an $n$-vertex, $(2k-1)$-partite graph $G$ has $o(n^{2\ell+1})$ copies of $C_{2\ell+1}$'s, then
it can be made $C_{2k-1}$-free by deleting $o(n^2)$ edges. 

\section{Algorithmic aspects}

The algorithmic version of Szemerédi's Regularity lemma proved to be very useful in many problems in computer science, e.g., by providing good approximation algorithms for several NP-complete questions. 
The graph functionals we considered can be used to show the existence of a large super-regular pair, but they are not capable of finding one effectively.
In this section we present a deterministic polynomial time algorithm for finding a large $\ep$-regular subgraph in a 
sufficiently dense graph. This 
algorithm can then be used for decomposing the edge set of a graph into large $\ep$-regulars pairs, similarly to Theorem~\ref{dekomp}, albeit the result will be somewhat weaker. Our algorithm is based on the method of~\cite{ADLRY}.


Let $M(n)$ denote the time needed to multiply two $n\times n$ matrices with $0, 1$ entries over the integers (so $M(n)=O(n^{2.376})$).   
For proving an algorithmic version of the Regularity lemma, the authors of~\cite{ADLRY}, among other lemmas, used the following\footnote{This is Corollary 3.3 in~\cite{ADLRY}, slightly rewritten.}:
 
\begin{lemma}\label{algi}
Let $H$ be a bipartite graph with equal parts $|A|=|B|=n.$ Let $2n^{-1/4}<\ep < 1/16.$ Then there is a $O(M(n))$ time deterministic algorithm that verifies that $H$ is $\ep$-regular, or finds
two subsets, $A_1\subset A,$ $B_1\subset B,$ $|A_1|, |B_1|\ge \ep^4 n/16,$ such that $|d(A, B)-d(A_1, B_1)|\ge \ep^4.$  
The algorithm can be parallelized and implemented in $NC^1.$  
\end{lemma}  
 
\medskip

We call the sets $A_1$ and $B_1$ the {\it witnesses of $\ep^4$-irregularity.} Note that the cardinalities of the witnesses for irregularity could  be much smaller than $\ep n.$ As it is proved in~\cite{ADLRY}, this is unavoidable unless $P=NP.$ 
We need a simple lemma before discussing the decomposition algorithm. 

\begin{lemma}\label{suru1}
Let $H$ be a bipartite graph with vertex parts $A, B$ such that $|A|=|B|=n,$ and let $0<\eta<d \le 1.$ Assume that the density of $H$ is $d,$ and that there exists $A_1\subset A, B_1\subset B$ 
such that $|A_1|, |B_1| \ge \eta n,$ and $|d-d(A_1, B_1)|\ge \eta.$ Then in polynomial time we can find $A'\subset A, B'\subset B$ with $|A'|=|B'|\ge \eta n$ such that $d(A', B')\ge d+\eta^3.$   
\end{lemma}

\noindent {\bf Proof:} We begin with the case $d(A_1, B_1)\le d-\eta.$ Let us assume first, that $|A_1|<n.$ Define the sets $A_2=A-A_1$ and $B_2=B-B_1.$ Below we show that one of the densities $d(A_1, B_2),$ $d(A_2, B_1)$ or $d(A_2, B_2)$
must be at least $d+\eta^3.$ 

Suppose not. Then $$e(H)=d n^2< (d-\eta) |A_1|\cdot |B_1|+(d+\eta^3)(n^2-|A_1|\cdot |B_1|).$$ This implies that $$\eta |A_1|\cdot |B_1| < \eta^3 n^2-\eta^3 |A_1|\cdot |B_1|.$$
Using that $|A_1|, |B_1|\ge \eta n,$ we arrived at a contradiction. Hence, one of the subgraphs $H[A_1, B_2],$ $H[A_2, B_1]$ or $H[A_2, B_2]$ must have density at least $d+\eta^3.$
Denote the parts of the densest of these subgraphs by $A''$ and $B'',$ here $A''\subset A, \ B''\subset B.$ 

Clearly, $|A''|, |B''|\ge \eta n.$ If $|A''|=|B''|,$ then we are done. Assume, that $|A''|>|B''|.$ Order the vertices of $A''$ in decreasing order according to their degrees in $B''.$ Keep the first $|B''|$
vertices, and discard the rest. Call the resulting set $A',$ and let $B'=B''.$ It is easy to see that $d(A'', B'')\le d(A', B'),$ so we are done.

There is one case left, when $|A_1|=n.$ Clearly, $|B_1|<n,$ so $B_2=B-B_1$ is non-empty. Assume that $d(A, B_2)< d+\eta^2.$ Then we would have 
$$d n^2 < (d-\eta) |A|\cdot |B_1| + (d+\eta^2) |A|\cdot |B_2|,$$ which implies that
$$ \eta |B_1|\cdot n < \eta^2 n^2 -\eta^2n|B_1|.$$ We arrived at a contradiction again, since $|B_1|\ge \eta n.$ By discarding those $|B_2|$ vertices of $A$ that have the smallest degrees into $B_2$
we get two equal sized subsets $A'$ and $B'=B_2$ such that $d(A', B')\ge d+\eta.$ 

Let us now assume that $d(A_1, B_1)\ge d+\eta.$ We are almost done. The only problem is if $|A_1|\neq |B_1|.$ Then we discard vertices from the larger set, similarly to the way it is described above, and we
obtain the pair of sets $A', B'$ such that $d(A', B')\ge d+\eta$ and $|A'|=|B'|\ge \eta n.$

It is easy to see that the above $A', \ B'$ subsets can be found in $O(n^2)$ time: we need to compute densities of at most four subgraphs, and order vertices according to their degrees. 
Hence, we proved what was desired.  
\qedf

\medskip

Using Lemma~\ref{suru1} we can easily formulate a polynomial time algorithm that finds a large $\ep$-regular subgraph in a graph. 
Let $G=(A, B; E)$ be a bipartite graph with $|A|=|B|=n$ and density $d,$ and assume that $0<\ep \ll d\ll 1,$ such that 
$\ep > exp({\frac{4\cdot 16^3 \log {1/\ep}}{\ep^{12}}}) 2n^{-1/4}.$ 

\begin{enumerate}

\item Apply Lemma~\ref{algi}. If $G[A, B]$ is $\ep$-regular, stop. 

\item If not, by Lemma~\ref{suru1}, substituting $\ep^4/16$ for $\eta,$ we can find a balanced subgraph $G[A', B']$ with $|A'|=|B'|\ge \ep^{12} |A|/16^3$ having density at least $d(G[A, B])+\ep^{12}/16^3.$ 

\item Let $A=A', \ B=B',$ and continue with Step 1.

\end{enumerate}

\smallskip

The above algorithm stops in at most $16^3/\ep^{12}$ steps, since if the density of a bipartite graph is 1, it must be $\ep$-regular.
With the above we have proved the following.

\begin{prop}\label{algProp}
Let $G=(V, E)$ be a balanced bipartite graph on $n$ vertices with density $d_G$ and $0<\varepsilon<1/4$ be a real number such that the following is satisfied:
$e^{r/4} 2n^{-1/4}<\ep \ll d_G\le 1,$ where 
$r= (16^3/\varepsilon^{12})\log (16^3/\varepsilon^{12}).$ Then $G$ contains an 
$\varepsilon$-regular subgraph $H\in \cB$
with density $\delta \ge d_G$ and $v(H)\ge e^{-r}n.$ 
\end{prop}

Iterating the above procedure as in the proof of Theorem~\ref{dekomp} we get the following.

\begin{theorem}\label{algdekomp}
Let $G$ be a balanced bipartite graph on $n$ vertices with density $d_G,$ and let $0<\delta \le 1$ and $0< \ep <1/4$ 
be real numbers such that the following is satisfied:
$e^{r/4} 2n^{-1/4}<\ep \ll \delta \le 1,$ where 
$r= (16^3/\varepsilon^{12})\log (16^3/\varepsilon^{12}).$
Then the \emph{edge set} of $G$ can be decomposed as follows: 
$E(G)$ can be written as the edge-disjoint union of the $\varepsilon$-regular balanced
bipartite graphs $F_1, \ldots, F_{K},$ and another balanced bipartite graph $F_0,$ where 
$K=K(\varepsilon, \delta, d_G)\le d_G e^r/\delta.$ For $i\ge 1$ each $F_i$ has at least 
$m=m(\varepsilon)\ge e^{-r}n$ vertices and density at least $\delta,$ 
while $F_0$ has density less than $\delta.$ 
\end{theorem}

\medskip

Note that in Proposition~\ref{algProp} we guarentee only $\ep$-regularity, not $(\ep, \delta)$-super-regularity for
some $\delta.$ If one needs super-regularity, Claim~\ref{szup-atlagfok} can be used, resulting only a small loss in
the size of the pair. Of course, this remark also
applies for Theorem~\ref{algdekomp}.
 
%

\end{document}